\documentclass[12pt]{article}
\usepackage[dvips]{color} 
\usepackage{amsmath, amsthm, amssymb}
\usepackage[left=1.5in,top=1in,right=1in,nohead]{geometry}
\usepackage{hyperref}
\usepackage{latexsym}
\usepackage{mathrsfs}
\usepackage{rawfonts}
\usepackage[dvips,dvipdfmx]{graphicx}
\usepackage{float}
\usepackage[all]{xy}
\usepackage[sc,osf]{mathpazo}  %\usepackage{kpfonts}

\usepackage{amsmath,amsthm,amssymb}
\usepackage{hyperref}
\usepackage{latexsym}
\usepackage{bmpsize}

\def\bct{\begin{center}}
\def\ect{\end{center}}
\def\bit{ \begin{itemize} }
\def\eit{ \end{itemize} }
\def\beg{\begin}

\def\<{\langle}
\def\>{\rangle}

\def\mbb{\mathbb}

\def\mbbr{\mathbb R}
\def\mbbz{\mathbb Z}

\def\ni{\noindent}

\def\tcb{\textcolor{blue}}

\def\tn{\textnormal}
\def\wt{\widetilde}

\newtheorem{thm}{Theorem}[section]

\newtheorem{theorem}{Theorem}
\theoremstyle{plain}

\numberwithin{equation}{section}

\title%[On some special submanifolds of Page space]
{On special
submanifolds of the Page space}
\author{Mustafa Kalafat \and Ramazan Sar\i}
\begin{document}
\maketitle

\begin{abstract}
In this paper, we study some classes of submanifolds of codimension one and two %surfaces 
in the Page space. 
These submanifolds are totally geodesic. 
%Among them, we illustrate that there are some totally geodesic examples ones and some minimal ones.  
We also compute their curvature and show that some of them are constant curvature spaces.  
Finally we give information on how the Page 
space is related to some other metrics on the same underlying smooth manifold. 
%These surfaces are examine totally geodesic manifolds with ?--> Gauss-Bonnet theorem. 
%After we are investigate minimal submanifolds. 
%These surfaces are calculated scalar curvature. 
%Which ones are ricci flat.
\end{abstract}

\maketitle

\section{Introduction}

\bigskip Let $M$ be a complex manifold with complex structure $J$. We say
that a Riemannian metric $g$ on a complex manifold $M$ is {\em Hermitian} if 
%$g(X,Y)=g(JX,JY)$ for any $X,Y\in T_{x}M$ for all $x\in M.$ 
$|X|_g=|JX|_g$ for any $X\in T_{x}M$ and $x\in M.$ 
Then the triplet $(M,J,g)$ is
called {\em Hermitian manifold.}  
A Riemannian metric is called {\em Einstein} if its Ricci tensor $\tn{Ric}$ is a constant multiple of itself, %the metric tensor, 
i.e. $\tn{Ric}=\lambda g$ for same constant $\lambda \in 
\mathbb{R}$ called the { Einstein constant.}  
If the metric is Einstein on the Hermitian manifold $M$, then we call it 
{\em Einstein-Hermitian.} 
%\vspace{.05in}
In %1978 
\cite{P:1978}, D. Page introduced an %studied a metric on 
Einstein-Hermitian metric on the blow up of the complex projective plane. %$4$-manifold. 
%This metric was called {\em Page metric.} Later, Page manifold called
%Einstein-Hermitian 4-Manifold.  
%M. Kalafat and C. Koca give mathematical exposition of the Page metric and introduce an
%efficient coordinate system for it \cite{ehpbisec}. 
%They show that Page metric does not have positive bisectional curvature.
%\vspace{.05in}
%\section{Page Metric} \bigskip 
%The Page metric was introduced by D. Page in 1978 
%It is introduced as a limiting metric of Kerr-de Sitter solution. 
The manifold is the nontrivial $S^{2}$\
fibre bundle over $S^{2}$\ and has $\chi =4$\ and $\tau =0.$ It is the
unique Einstein-Hermitian non K\"{a}hler metric on the blow up of complex
projective plane. One can first think of the following metric on the product $%
S^{3}\times I$ where $I$ is the closed interval $[0,\pi ]$

\begin{equation} \label{page1}
g_{Page}:=V(r)dr^{2}+f(r)(\sigma _{1}^{2}+\sigma _{2}^{2})+\frac{C\sin ^{2}r}{V(r)}
\sigma_3^2
\end{equation}

%$$g_{Page}:=V(r)dr^{2}+f(r)(\sigma _{1}^{2}+\sigma _{2}^{2})+\frac{C\sin ^{2}r}{V(r)}
%\sigma _{3}^{2}$$
\ni where the functions
\begin{eqnarray*}
V(r) &=&\frac{1-a^{2}\cos ^{2}r}{3-a^{2}-a^{2}(1+a^{2})\cos ^{2}r} \\
f(r) &=&4\,{ 1-a^{2}\cos ^{2}r \over 3+6a^{2}-a^{4}} \\
C &=&\left( \frac{2}{3+a^{2}}\right) ^{2}
\end{eqnarray*}%
$\sigma _{1},\sigma _{2},\sigma _{3}$ are standard left invariant 1-forms on
the Lie group $SU(2)\thickapprox S^{3}$ and $a$ is a small, positive constant 
which can approximately be computed as $a\approx 0.28170$ (So $C\approx 0.42183$). 
It is actually defined to be the unique positive root of the following quartic polynomial,
$$a^{4}+4a^{3}-6a^{2}+12a-3=0.$$

\begin{figure}[h] \bct 
\includegraphics[width=.7\textwidth]{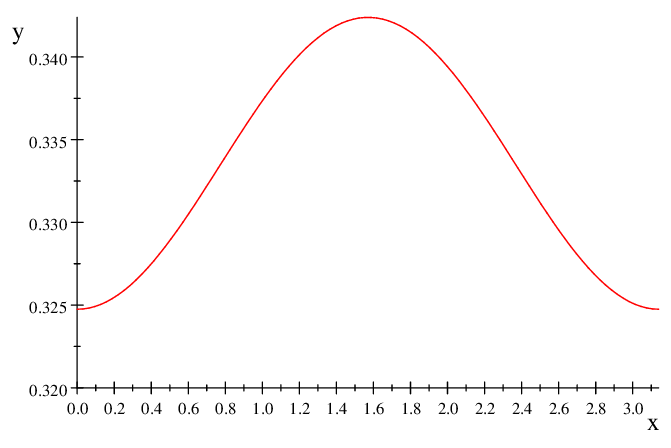}\\ 
\ect  \caption{\em The function $V(r)$ for $r\in [0,\pi]$. 
Minimum values are $V(0,\pi)\approx 0.324776$ and the maximum value $V(\pi/2)\approx 0.342397$. } \label{functionv}
\end{figure}

\begin{figure}[h] \bct 
\includegraphics[width=.7\textwidth]{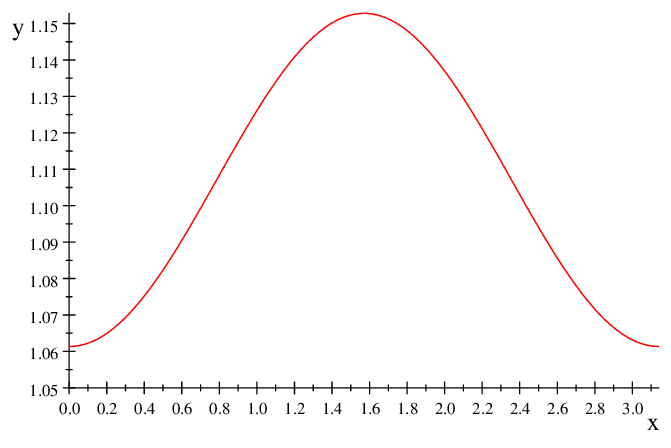}\\ 
\ect  \caption{\em The function $f(r)$ for $r\in [0,\pi]$. 
Here, at the endpoints $f(0)=f(\pi)\approx  1.061462$, 
the maximum value is approximately equal to $1.152811$, taken in the middle.}
\label{functionf} 
\end{figure}

\begin{figure}[h] \bct 
\includegraphics[width=.7\textwidth]{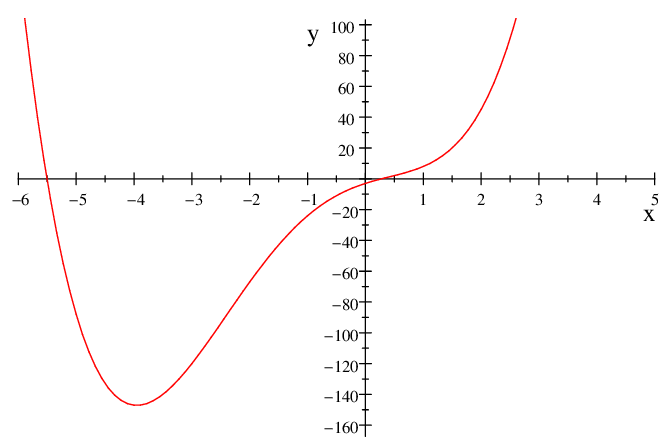}\\ 
\ect  \caption{\em The polynomial $p(x)=x^4+4x^3-6x^2+12x-3$. } \label{polynomialfora}
\end{figure}

\vspace{.05in}

At the end points $r=0$ and $r=\pi $, the metric shrinks to the the round
metric on $S^{2}.$ Thus, $g$ descends to the metric $g_{Page}$ on the
quotient $(S^{3}\times I)/\sim $ where $\sim $ identifies the fibers of the
Hopf fibration $h:S^{3}\rightarrow S^{2}$ on the two ends $S^{3}\times \{0\}$
and $S^{3}\times \{\pi \}$ of the cylinder $S^{3}\times I.$ The resulting
manifold is indeed the connected sum 
$\mbb{CP}_2\sharp\, \overline{\mbb{CP}}_2$.   
See \cite{ehpbisec} for a mathematical exposition of the Page metric for
further details. There, an efficient coordinate system is introduced and studied. 
This coordinate system uses Euler
angles %\cite{xiaozhang}
 on the $S^3\subset\mbbr^4$ which e.g. realize the Hopf fibration in the best.

\bigskip 
\begin{equation*}
0\leq \theta \leq \pi ,%
\begin{array}{cc}
& 
\end{array}%
0\leq \phi \leq 2\pi ,%
\begin{array}{cc}
& 
\end{array}%
0\leq \psi \leq 4\pi
\end{equation*}

\begin{eqnarray}
x_{1} &:=&r\cos \frac{\theta }{2}\cos \frac{\psi +\phi }{2} \nonumber\\
x_{2} &:=&r\cos \frac{\theta }{2}\sin \frac{\psi +\phi }{2} \nonumber\\
x_{3} &:=&r\sin \frac{\theta }{2}\cos \frac{\psi -\phi }{2}\label{eulercoordinates}\\
x_{0} &:=&r\sin \frac{\theta }{2}\sin \frac{\psi -\phi }{2} \nonumber
\end{eqnarray}
where the Hopf fibration in these coordinates is just a projection \cite{fip},%
\begin{equation*}
h:S^{3}\rightarrow S^{2},%
\begin{array}{cc}
& 
\end{array}%
h(\theta ,\psi ,\phi )=(-\phi ,\theta ).
\end{equation*}%
Here the exchange $\phi \leftrightarrow \theta $ is needed to relate to the
calculus angles on $S^{2}.$ Changing $\psi $ does not change the element in
the image. So, whenever the image $\phi ,\theta $ is fixed, $\psi $
parametrizes the Hopf circle.  
An %orthonormal, 
invariant coframe %$\{\sigma _{1},\sigma _{2},\sigma _{3}\}$
on $S^{3}$ is given as follows:%
\begin{eqnarray*}
\sigma _{1} &=&(\sin \psi \,d\theta -\sin \theta \cos \psi \,d\phi )/2 \\
\sigma _{2} &=&(-\cos \psi \,d\theta -\sin \theta \sin \psi \,d\phi )/2 \\
\sigma _{3} &=&(d\psi +\cos \theta \,d\phi )/2
\end{eqnarray*} 
%On can check in a straightforward manner the identities%
%\begin{equation*}
%d\sigma _{i}=2\sigma _{i+1}\wedge \sigma _{i+2} ~(mod ~3)
%\end{equation*}%
%and%
%\begin{equation*}
%\sigma _{1}^{2}+\sigma _{2}^{2}=(d\theta ^{2}+\sin ^{2}\theta \,d\phi ^{2})/4
%\end{equation*}%
In terms of these new coordinates %Plugging these into the Page 
the metric becomes, %'s expression we get

\bigskip

\begin{equation}\label{pageeuler}
\hspace{-16mm}g_{\mathrm{Page}}=Vdr^2+\left\{ {f\over 4}\sin^2\theta +{C\sin^2r \cos^2\theta\over 4V(r)} \right\}\hspace{-1mm}d\phi^2+ {C\sin^2r\over 4V(r)}d\psi^2
+ {C\sin^2r \cos\theta\over 4V(r)}
(d\psi\otimes d\phi+d\phi\otimes d\psi)+{f\over 4}d\theta^2.\end{equation}

%$$\hspace{-16mm}g_{\mathrm{Page}}=Vdr^2+\left\{ {f\over 4}\sin^2\theta +{C%\sin^2r \cos^2\theta\over 4V(r)} \right\}\hspace{-1mm}d\phi^2+ {C\sin^2r\over 4V(r)}d\psi^2
%+ {C\sin^2r \cos\theta\over 4V(r)}
%(d\psi\otimes d\phi+d\phi\otimes d\psi)+{f\over 4}d\theta^2.$$
%\begin{equation*}
%g_{page}=Vdr^{2}+\frac{C\sin ^{2}r}{4V(r)}d\psi ^{2}+\{\frac{f}{4}\sin
%^{2}\theta +\frac{C\sin ^{2}r\cos ^{2}\theta }{4V(r)}\}d\phi ^{2}+\frac{%
%C\sin ^{2}r\cos ^{2}\theta }{4V(r)}(d\psi \otimes d\phi +d\phi \otimes d\psi
%)+\frac{f}{4}d\theta ^{2}
%\end{equation*}%

\ni Letting $U:=\sqrt{V(r)},\,h:=\sqrt{f(r)},\, D:=\sqrt{C}$ we choose the simplest Vierbein 
i.e. orthonormal coframe as follows 

\begin{equation}\label{vierbein}
\{e^{0},e^{1},e^{2},e^{3}\}:=\{Udr,\, 2^{-1}h\sin\theta \,d\phi,\,2^{-1}U^{-1} D\sin r \,(d\psi+\cos\theta \,d\phi),\, 
2^{-1} h\,d\theta \}.\end{equation}

\ni After this, using the structure equations \cite{egh}\, $de^a=-\omega^a_{\phantom{1} b}\wedge e^b$ \, for \, $0\leq a\leq 3$, %non-trivial 
connection $1$-forms %$\omega^a_{\phantom{1} b}$ 
of the Page metric can be computed as the following. % \cite{ehpbisec}  
%\begin{equation}
%\omega^1_{\phantom{1} 0}=U^{-2}(U\cot r-\dot U)\,e^1,%
%\begin{array}{cc}& \end{array}
%\omega^2_{\phantom{1} 0}=\dot hh^{-1}U^{-1}\,e^2,
%\begin{array}{cc}& \end{array}
%\omega^3_{\phantom{1} 0}=\dot hh^{-1}U^{-1}\,e^3
%\end{equation}
%\begin{equation*} \omega^2_{\phantom{1} 3}=2h^{-1}\cot \theta\,e^2. \end{equation*}
%{ \renewcommand*{\arraystretch}{1.4} \hspace{-2.5cm}
%\begin{equation}\label{ambientconnectiononeforms}\hspace{-2.5cm}
%\begin{array}{ccl}
%\omega^1_{\phantom{1} 0}=U^{-1}h^{-1}\dot h \,e^1 + 2^{-1}U^{-1}h^{-1}\dot h \sin^{-1}%\theta\,e^3,~~ &
%\omega^2_{\phantom{1} 0}=\sin^{-1}r\,(U^{-1}\sin r)_r\,e^2,~~ &
%\omega^3_{\phantom{1} 0}=2^{-1}U^{-1}h^{-1}\dot h \sin^{-1}\theta\,e^1,\\ 
%\omega^2_{\phantom{1} 1}=DU^{-1}h^{-2}\sin r\,e^3, &
%\omega^3_{\phantom{1} 2}=2 D U^{-1} h^{-2} \sin r \,e^1, &
%\omega^1_{\phantom{1} 3}=2^{-1}U^{-1}h^{-1}\dot h \sin^{-1}\theta\,e^0+
%2h^{-1}\cot\theta\,e^1 - DU^{-1}h^{-2}\sin r\,e^.
%\end{array}
%\end{equation}   }

{ \renewcommand*{\arraystretch}{1.4} 
\begin{equation}\label{ambientconnectiononeforms}\hspace{-2.2cm}
\begin{array}{r}
\omega^1_{\phantom{1} 0}=U^{-1}h^{-1}\dot h \,e^1 + 2^{-1}U^{-1}h^{-1}\dot h \sin^{-1}\theta\,e^3,~
\omega^2_{\phantom{1} 0}=\sin^{-1}r\,(U^{-1}\sin r)_r\,e^2,~
\omega^3_{\phantom{1} 0}=2^{-1}U^{-1}h^{-1}\dot h \sin^{-1}\theta\,e^1,\\ 
\omega^2_{\phantom{1} 1}=DU^{-1}h^{-2}\sin r\,e^3,~~ 
\omega^1_{\phantom{1} 3}=2^{-1}U^{-1}h^{-1}\dot h \sin^{-1}\theta\,e^0+
2h^{-1}\cot\theta\,e^1 - DU^{-1}h^{-2}\sin r\,e^2,\\
\omega^3_{\phantom{1} 2}= D U^{-1} h^{-2} \sin r \,e^1.
\end{array}
\end{equation}
}

\ni See Appendix \ref{ode} for details. Note that these $1$-forms constitute a skew-symmetric matrix of forms. 
One can further go and compute the curvature $2$-forms using the identity 
$R^a_{\phantom{1} b}:=d\omega^a_{\phantom{1} b}+\omega^a_{\phantom{1} c}\wedge \omega^c_{\phantom{1} b}$. These make up a skew-symmetric matrix of $2$-forms as well. 
Looking at the coefficients with respect to a basis yields the Riemann curvature tensor coefficients 
$R^a_{\phantom{1} b}={1\over 2}R^a_{\phantom{1} bcd}\,e^c\wedge e^d.$

\vspace{3mm}

\ni We summarize some of our results as follows. 

\begin{thm} There are totally geodesic subsurfaces of the Page space which are 
 isometric to tori and sphere with their constant zero and positive  curvature metrics respectively.
\end{thm}

\begin{thm} There are totally geodesic hypersurfaces of the Page space which are 
diffeomorphic to $S^3$ and $S^1\times S^2$.
\end{thm}

\vspace{.05in}

%In this paper, we are studied submanifolds of Page manifold. Firstly we
%introduce the concept of Page metric. After we are given an efficient Euler
%coordinate system for Page metric. In section 3, we are determined
%submanifolds of Page manifold.\ These surfaces are examine totally
%geodesic submanifolds with Gauss-Bonnet theorem.In section 4, we are
%investigate minimal submanifolds of Page space. Last section, submanifolds
%are calculated curvatures. We are proved \ Ricci flat some of submanifold.

In sections $\S$\ref{secsurfaces} and $\S$\ref{sechypersurfaces} 
we study the surfaces and hypersurfaces inside the Page space. In $\S$\ref{secyamabe} 
we compare it with the other metrics and study the relation with the Yamabe problem. 
Finally in Appendix $\S$\ref{ode} we derive the connection $1$-forms. 
  
\vspace{.05in}

\noindent{\bf Acknowledgements.} Thanks to C. Koca and \" O. Kelek\c ci for useful comments. 
This work is partially supported by T\"ubitak (Turkish science and research council) grant {$\sharp$}113F159. The Ansi C programing language \cite{kernighantritchie} 
is used for the computation of some min/max results. Thanks to Alphan Es 
for helping with the graphics.

%\newpage

\section{Minimal surfaces inside the Page Space and their curvature}\label{secsurfaces}

Now we are in a position to analyse some of the submanifolds with partially constant Euler coordinates. Here, keep in mind that $\psi$ is the Hopf circle direction. 

\vspace{.05in}

$\tcb{\mathbf{S_1:}}$ By keeping $r\in(0,\pi),\,\theta\in[0,\pi] $ fixed and varying $\psi ,\phi $ one
obtains {\em tori}. At the ends for $r=0,\pi$, this construction produces {\em circles}.  
The restriction of the metric to these tori is computed as
$$\hspace{-4mm}g_{S_1}=\{\frac{f}{4}\sin ^{2}\theta +\frac{C\sin ^{2}r\cos ^{2}\theta }{%
4V(r)}\}d\phi ^{2}+\frac{C\sin ^{2}r}{4V(r)}d\psi ^{2} 
+\frac{C\sin ^{2}r\cos ^{2}\theta }{4V(r)}(d\psi \otimes d\phi +d\phi
\otimes d\psi).$$
%$$g_{page} =\{\frac{f}{4}\sin ^{2}\theta +\frac{C\sin ^{2}r\cos ^{2}\theta }{%
%4V(r)}\}d\phi ^{2}+\frac{C\sin ^{2}r}{4V(r)}d\psi ^{2} 
%+\frac{C\sin ^{2}r\cos ^{2}\theta }{4V(r)}(d\psi \otimes d\phi +d\phi
%\otimes d\psi ).$$
We can work with the orthonormal coframe 
$\{e^1,e^2\}$ %=\{ \frac{h}{2}\sin \theta d\phi , \frac{D\sin r}{2U}(d\psi +\cos \theta d\phi )\}.
selecting from (\ref{vierbein}) on these surfaces since this spans the $\phi\psi$ plane. 
%So %On the other hand, %\tcr{bu hangi baza g\"ore ?:} 
%\begin{equation*}
%g_{ij}=\left( 
%\begin{array}{cc}
%\frac{f}{4}\sin ^{2}\theta +\frac{C\sin ^{2}r\cos ^{2}\theta }{4V(r)} & 
%\frac{C\sin ^{2}r\cos ^{2}\theta }{4V(r)} \\ 
%\frac{C\sin ^{2}r\cos ^{2}\theta }{4V(r)} & \frac{C\sin ^{2}r}{4V(r)}%
%\end{array}%
%\right) 
%\end{equation*}% then 
%we compute the volume form, % as, 
%\begin{equation*}
%\left\vert g_{ij}\right\vert =\frac{f.C}{16V}\sin ^{2}\theta \sin ^{2}r
%\end{equation*}%
%\begin{equation*}dV=e^1\wedge e^2=\frac{Dh\sin \theta \sin r}{4U}\,d\psi \wedge d\phi\end{equation*}
Connection $1$-form can be computed from the structure equations 
$0=de^i+\wt{\omega}^i_{\phantom{1} j}\wedge e^j$. Both components vanish hence 
$\wt{\omega}^1_{\phantom{1} 2}=0$.

$$\wt R^1_{\phantom{1}2} =d\wt{\omega}_{\phantom{1} 2}^{1}+\wt{\omega}_{\phantom{1} i}^{1}\wedge \wt{\omega}_{\phantom{1} 2}^i=0.$$
%+w_{1}^{1}\wedge
%w_{2}^{1}+w_{2}^{1}\wedge w_{2}^{2}+w_{3}^{1}\wedge w_{2}^{3} \\
%&=&(\sqrt{V}\cot r\frac{\sqrt{C}\sin r}{2\sqrt{V}}d\psi )\wedge (\frac{dr}{2%
%\sqrt{f}}\frac{1}{\sqrt{f}}\frac{1}{\sqrt{V}}\frac{\sqrt{f}}{2}\sin \theta
%(d\psi +d\phi )) \\
%&=&0
%\end{eqnarray*}%
$$\wt R_{\phantom{1}2}^{1}=\frac{1}{2}\,\wt R_{\phantom{1}212}^{1}\,e^{1}\wedge e^{2}
~~~\tn{implies that}~~~ %\Rightarrow~ 
\wt R_{1212}=\wt R_{\phantom{1}212}^{1}=0.$$
%\begin{equation*}\int_{S_1} R_{212}^{1}\,dV=0.\end{equation*}
\ni So that all these tori are {\em flat}. We know that the ambient connection 
$\omega^1_{\phantom{1} 2}|_{S_1}=0$ so that the second fundamental form, 
which is the difference of two connections, vanishes on these 
subtori. Hence these tori are {\em totally geodesic}, consequently minimal.   
One can alternatively compute Christoffel symbols from metric coefficients 
directly to see that they all vanish to yield zero curvature.  
The metric is conformally equivalent to
%$$\wt g_{S_1}=\{fV\sin ^2\theta +C\sin^2 r\cos^2\theta\}d\phi^2
%+C\sin^2 r d\psi^2 
%+C\sin^2r\cos^2\theta (d\psi \otimes d\phi +d\phi\otimes d\psi).$$
%$$\wt g_{S_1}=\{fV\sin^{-2}r \sin^2\theta +C\cos^2\theta\}d\phi^2
%+C d\psi^2 + C\cos^2\theta (d\psi \otimes d\phi +d\phi \otimes d\psi).$$
$$\wt g_{S_1}=\{C^{-1}fV\sin^{-2}r \sin^2\theta+\cos^2\theta\}d\phi^2
+ d\psi^2 + \cos^2\theta (d\psi \otimes d\phi +d\phi \otimes d\psi).$$
%$$\wt g_{ij}=\left[ 
%\begin{array}{cc}
%C^{-1}fV\sin^{-2}r \sin^2\theta + \cos^2\theta & \cos^2\theta\\
%\cos^2\theta & 1\end{array}\right] $$
$$\wt g_{ij}=\left[ 
\begin{array}{cc}
C^{-1}fV\sin^{-2}r \tan^2\theta + 1 & 1\\
1 & \cos^{-2}\theta\end{array}\right] $$

\ni We can compute the angle of the parallelogram as follows.

$$\cos\Theta=g_{\phi\psi}/\sqrt{g_{\phi\phi}g_{\psi\psi}}
=\cos\theta/\sqrt{C^{-1}fV\sin^{-2}r \tan^2\theta + 1}$$

\ni The radius function is the following,

%$$L=\int_0^{4\pi}|\partial_\psi|d\psi/\int_0^{2\pi}|\partial_\phi|d\phi=
%2g_{\psi\psi}/g_{\phi\phi}=$$
$$\tn{R}=\int_0^{2\pi}|\partial_\phi|\,d\phi\, /\hspace{-1mm}\int_0^{4\pi}|\partial_\psi|\,d\psi
=2^{-1}g_{\phi\phi}/g_{\psi\psi}
=2^{-1}\{C^{-1}fV\sin^{-2}r \sin^2\theta + \cos^2\theta\}.$$

\ni If we compute he minimum and maximum of the radius and cosine 
while varying $r\in(0,\pi),\theta\in[0,\pi]$, we got the following intervals,  
$$R \in [0.408520,\infty) ~~\tn{and}~~ \cos\Theta \in [-1,1].$$

\begin{figure}[h] \bct 
\includegraphics[scale=0.7]{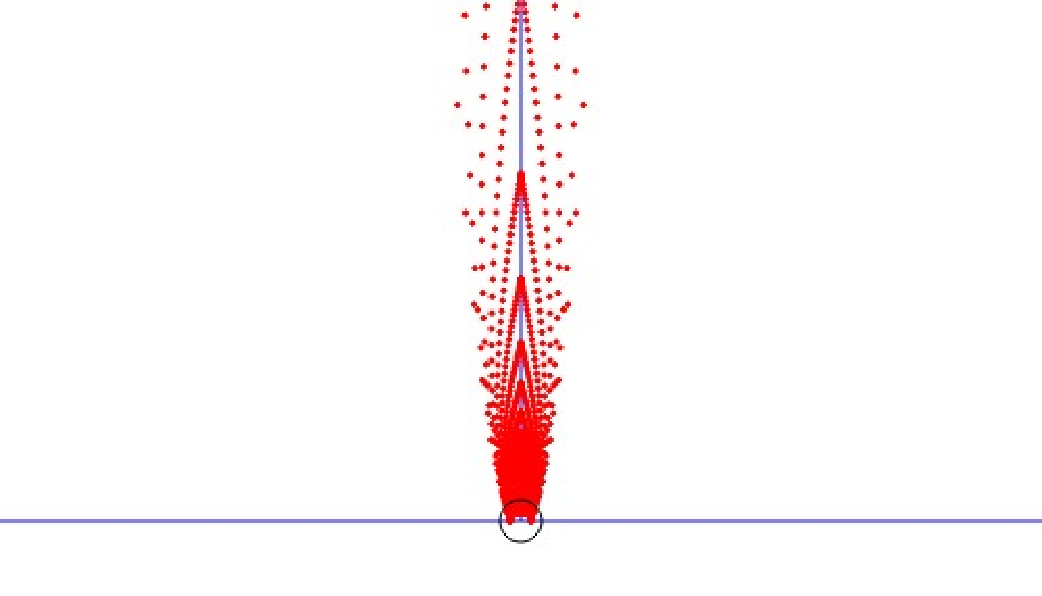}
%,natwidth=200,natheight=200   also bb=0 0 1280 960
\ect  \caption{\em Image of the $\psi,\phi$ tori $S_1$ on the Teichm\" uller %and Moduli 
space.}  \label{teichmullermodulis1}
\end{figure}

\ni where the minimum radius is taken at $r=1.571425, \theta=1.571425$. 
See Figure \ref{teichmullermodulis1} for the coverage of these metrics on the Teichm\" uller 
space. Note that the whole moduli space is covered with these classes of flat metrics.   

%MINRADIUS=0.408520 at r=1.571425 teta=1.571425, MAXRADIUS=17917768892416.000000 at r=3.141593 teta=1.570796
%maximumu sonsuza gidiyo yani

%MINCOSTETA=-1.000000 at r=3.140964 teta=3.141593, MAXCOSTETA=1.000000 at r=1.970407 teta=0.001257
%da buldu benim C

%\begin{figure}[h] \bct 
%\advance\leftskip-3cm
%\advance\rightskip-3cm
%%\includegraphics[width=.7\textwidth]{pageteichmullers200-200-x10.eps}\\ 
%%\includegraphics[bb=0 0 160 160,scale=0.7]{pageteichmullers200-200-x10.eps}
%\includegraphics[bb=-260 -260 160 160,scale=0.7]{pageteichmullers200-200-%x10.eps}
%%,natwidth=200,natheight=200   also bb=0 0 1280 960
%\ect  \caption{\em Image of the $\psi,\phi$ tori $S_1$ on the Teichm\" uller %%and Moduli 
%space.}  \label{teichmullermodulis1}
%\end{figure}

%\cite{jostriemannsurface} 

$\tcb{\mathbf{S_2:}}$  This time keep $\psi ,\phi $ fixed and vary $r,\theta $ to
obtain a {\em cylinder}. Its two ends lie on the two Hopf spheres.      %n embedded rectangle. %\tcr{what ?}. 
Then, plugging these into the Page metric's expression we get%
\begin{equation*}
g_{S_2}=Vdr^2+\frac{f}{4}d\theta ^{2}
\end{equation*}%
and using the orthonormal coframe
$\{e^0,e^3\}$ %=\{Udr,\frac{h}{2}d\theta \}.
%\begin{equation*}
%\left\vert g_{ij}\right\vert =\frac{V.f}{4}
%\end{equation*}%
%\begin{equation*}dV=\frac{\sqrt{V.f}}{2}\,dr\wedge d\theta. \end{equation*}
from the two structure equations we get,

{ \renewcommand*{\arraystretch}{1.4}
$$\begin{array}{l}
0=de^0+\wt\omega^0_{\phantom{1}3} \wedge e^3= 0 + \wt\omega^0_{\phantom{1}3} \wedge e^3
%\wt\omega^0_{\phantom{1}3} \wedge 2^{-1}hd\theta
\\
0=de^3+\wt\omega^3_{\phantom{1}0} \wedge e^0=2^{-1}\dot h \,dr\wedge d\theta + 
\wt\omega^3_{\phantom{1}0} \wedge Udr
\end{array}$$ }

\ni we conclude that ~~ 
$\wt\omega^3_{\phantom{1}0}=2^{-1}U^{-1}h\,d\theta=U^{-1}h^{-1}\dot h\, e^3$, 
~ which is the same as the one in the ambient space so that  these cylinders are 
{\em totally geodesic.} Computing the curvature $2$-form, 
{ \renewcommand*{\arraystretch}{1.4}
\begin{eqnarray*}
\wt R^0_{\phantom{1}3} &=&d\wt\omega^0_{\phantom{1}3}+\wt\omega^0_{\phantom{1}i}\wedge \wt\omega^i_{\phantom{1}3}\\
&=&2^{-1}(U^{-2}\dot U\dot h-U^{-1}\dot h) \,dr\wedge d\theta \\
&=& (U^{-3}\dot U-U^{-2}h^{-1}\dot h)  \,e^0\wedge e^3.
\end{eqnarray*}
} 
\ni So that the Gaussian curvature is equal to ~$\wt R_{0303}=\wt R^0_{\phantom{1}303}=U^{-3}h^{-1}\dot U\dot h-U^{-2}h^{-1}\dot h$.

%\begin{equation*}
%R_{3}^{0}=\frac{1}{2}\,R_{303}^{0}\,e^{0}\wedge e^{3}\Rightarrow R_{303}^{0}=0
%\end{equation*}%
%\begin{equation*} \int_{S_2} R_{303}^{0}\,dV=0 \end{equation*}

\vspace{.05in}

\begin{figure}[h] \bct 
\includegraphics[width=.7\textwidth]{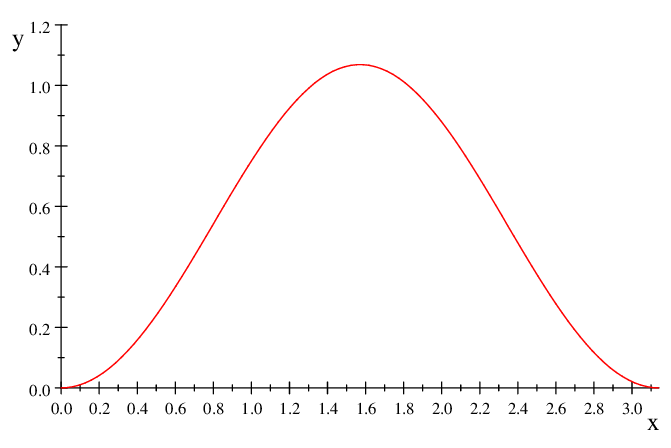}\\ 
\ect  \caption{\em The coefficient function    
${C \sin^2 r / fV}$​ of the conformally rescaled metric on the torus for $r\in [0,\pi]$.  This function vanishes at the endpoints and takes the approximate maximum value $1.068802$ at the midpoint.}  \label{konformalfaktors3}
\end{figure}

\begin{figure}[h!] \bct 
\includegraphics[width=.7\textwidth]{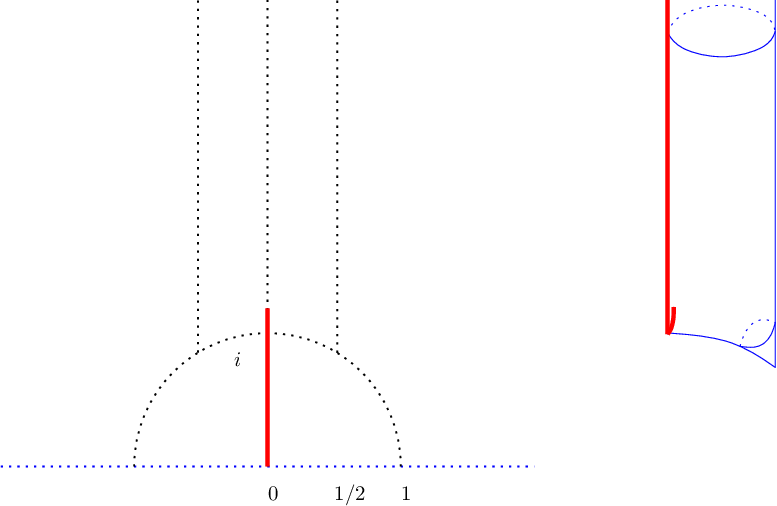}\\ 
\ect  \caption{\em Image of the $\psi,\theta$ tori $S_3$ on the Teichm\" uller and Moduli space.}  \label{teichmullermodulis3}
\end{figure}

$\tcb{\mathbf{S_3:}}$ We keep $r,\phi $ fixed and vary $\psi ,\theta $ to 
obtain an embedded cylinder. If we pass beyond the chart and vary 
$\theta\in [0,2\pi]$, then we get a {\em torus}.  At the endpoints for $r=0,\pi$ 
these tori degenerate to a circle parametrized by the angle $\theta$.   
Metric on the tori is a product, so that these tori are flat. % reads
%$$g_{S_3}=\frac{C\sin ^{2}r}{4V(r)}d\psi ^{2}+\frac{f}{4}d\theta ^{2}.$$
%and orthonormal frame$\ $
%\begin{equation*}
%\{e^{2},e^{3}\}=\{\frac{D\sin r}{2U}d\psi ,\frac{h}{2}d\theta \}
%\end{equation*}
%\begin{equation*}
%\left\vert g_{ij}\right\vert =\frac{Cf\sin ^{2}r}{16V}
%\end{equation*}%
%\begin{equation*}
%dV=\frac{1}{4}\sqrt{\frac{Cf}{V}}\sin r \, d\psi\wedge d\theta. 
%\end{equation*}
\ni %From the two structure equations $\wt\omega^2_{\phantom{1}3}=0$, 
%This is a product metric, so that these tori are flat. 
% ...........but not totally geodesic......... 
If we conformally rescale its metric we get, 
$$\wt g_{S_3}={C\sin^2 r \over f(r)V(r)}d\psi ^{2} + d\theta ^{2}.$$
The graph of the  coefficient function can be seen in Figure \ref{konformalfaktors3}. 
The range of this function is approximately in the open interval $(0, 1.068802)$ for $r\in (0,\pi)$. So that if one looks at the image of these tori on the Teichm\" uller space  
\cite{jostriemannsurface} that is a line on the imaginary axis of height slightly above $i$. The part under $i$ can be carried onto the moduli space via the linear fractional  transformation $f(z)=-1/z \in \tn{PSL}(2,\mbbz)$. So that it covers the whole imaginary axis above $i$. See Figure \ref{teichmullermodulis3}. 

\vspace{.05in}

$\tcb{\mathbf{S_4:}}$ We keep $\theta ,\psi $ fixed and vary $r,\phi $ to
obtain a {\em cylinder} with boundary components on the Hopf spheres at the two ends. 
Metric on this cylinder is 
\begin{equation*}
g_{S_4}=Vdr^{2}+\{\frac{f}{4}\sin ^{2}\theta +\frac{C\sin ^{2}r\cos
^{2}\theta }{4V(r)}\}d\phi ^{2}
\end{equation*}%
and orthonormal frame 
$\{e^{0},e^{1}\}$, %=\{Udr,\frac{h}{2}\sin \theta \,d\phi\}
from the two structure equations,

{ \renewcommand*{\arraystretch}{1.4}
$$\begin{array}{l}
0=de^0+\wt\omega^0_{\phantom{1}1} \wedge e^1=\wt\omega^0_{\phantom{1}1} \wedge e^1 
%2^{-1}hd\theta
\\
0=de^1+\wt\omega^1_{\phantom{1}0} \wedge e^0
=2^{-1}\dot h \sin\theta \,dr\wedge d\phi+\wt\omega^1_{\phantom{1}0} \wedge Udr
\end{array}$$ }

\ni we conclude that ~~ 
$\wt\omega^1_{\phantom{1}0}=1^{-1}U^{-1}\dot h \sin\theta\,d\phi=U^{-1}h^{-1}\dot h\, 
e^1$, ~ which is the same as the one in the ambient space so that  these cylinders are 
{\em totally geodesic.} Computing the curvature $2$-form, 
{ \renewcommand*{\arraystretch}{1.4}
\begin{eqnarray*}
\wt R^0_{\phantom{1}1} &=&d\wt\omega^0_{\phantom{1}1}+\wt\omega^0_{\phantom{1}i}\wedge 
\wt\omega^i_{\phantom{1}1}\\
&=&2^{-1}(U^{-2}\dot Uh-U^{-1}\dot h) \,dr\wedge d\theta \\
&=& (U^{-3}\dot U-U^{-2}h^{-1}\dot h)  \,e^0\wedge e^1.
\end{eqnarray*}
}
So that the Gaussian curvature is equal to ~$\wt R_{0101}=\wt R^0_{\phantom{1}101}=U^{-3}
\dot U-U^{-2}h^{-1}\dot h$.

\vspace{.1in}

$\tcb{\mathbf{S_5:}}$  Now, keep $\phi ,\theta $ fixed, vary $r,\psi $ to obtain
{\em spheres}. 
\begin{figure}[h!] \bct
\includegraphics[width=.4\textwidth]{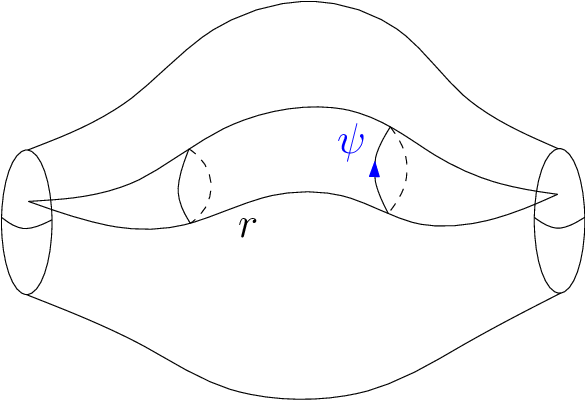}\\ \ect
  \caption{\em Varying $r$ and $\psi$.}  \label{sphere1}
\end{figure}
Metric restricted to these spheres is computed as 
\begin{equation*}
g_{S_5}=Vdr^{2}+\frac{C\sin ^{2}r}{4V(r)}d\psi ^{2}
\end{equation*}%
and using the orthonormal coframe 
$$\{e^{0},e^{2}\}=\{Udr,\frac{D\sin r}{2U}d\psi \},$$

%on the other hand%
%\begin{equation*}\left\vert g_{ij}\right\vert =\frac{C\sin ^{2}r}{4}\end{equation*}%
%then \begin{equation*}dV={D \over 2}\,\sin r\,dr\wedge d\psi .\end{equation*}

\ni we compute the connection $1$-form,
$$\wt\omega_{\phantom{1}0}^2=\cos\theta\, d\phi
=U^{-2}(-\dot U+U \cot r)\,e^2$$
which is the same as the ambient one so that the second fundamental form vanishes and hence these spheres are {\em totally geodesic} and {\em complex}, see \cite{ehpbisec}. The curvature,   
{ \renewcommand*{\arraystretch}{1.4}
$$\begin{array}{ccl}
\wt R^{2}_{\phantom{1}0} &=&d\wt\omega_{\phantom{1}0}^{2}+\wt\omega_{\phantom{1}i}^{2}
\wedge \wt\omega^i_{\phantom{1}0}\\
&=& 2^{-1}DU^{-3}\{3\dot U(-U^{-1}\dot U\sin r+\cos r)+(\ddot U\sin r+U\sin r)
  \}\,dr \wedge d\psi .
\end{array}$$
}
\begin{equation*}
\wt R^{2}_{\phantom{1}0}=\frac{1}{2}\, \wt R^{2}_{\phantom{1}020}\,e^2\wedge e^0
\end{equation*}
Then,
\begin{equation*}
\wt R^2_{\phantom{1}020}=U^{-3}\{3\dot U(-U^{-1}\dot U  +\cot r)+   \ddot U +U
  \}\, , %\cot \theta 
\end{equation*}
\ni and one can check through Gauss-Bonnet that, 
$$\int_{S_5} R_{020}^{2}\,dV
=\int_{S_5} 2^{-1}D\sin r\,  U^{-3}\{3\dot U(-U^{-1}\dot U  +\cot r)+   \ddot U +U \}\, dr\wedge d\psi 
=4\pi.$$

\vspace{.05in}

$\tcb{\mathbf{S_6:}}$  We keep $r,\psi $ fixed and vary $\phi ,\theta $ to obtain
{\em spheres}, again. Then restricting the Page metric we get 
\begin{equation*}
g_{S_6}=\{\frac{f}{4}\sin ^{2}\theta +\frac{C\sin ^{2}r\cos ^{2}\theta }{%
4V(r)}\}d\phi ^{2}+\frac{f}{4}d\theta ^{2}
\end{equation*}%
and working with the orthonormal coframe
$\{e^{1},e^{3}\}$ %=\{\frac{h}{2}\sin \theta\, d\phi,\,\frac{h}{2}d\theta\}.
%On the other hand 
%\begin{equation*}dV={f\sin\theta \over 4}\,d\phi\wedge d\theta % ~~~~\tcr{(form~ eklendi)}\end{equation*}%
we compute the connection $1$-form,
$$\wt\omega^1_{\phantom{1}3}=\cos\theta\, d\phi= 2h^{-1}\cot\theta\,e^1$$
is the same as the ambient one so that the second fundamental form vanishes and hence 
these spheres are totally geodesic. The curvature,   
\begin{eqnarray*}
\wt R_{\phantom{1}3}^{1} &=&d\wt\omega_{\phantom{1}3}^{1}+\wt\omega_{\phantom{1}i}^{1}
\wedge \wt\omega^i_{\phantom{1}3}\\
&=&-\sin \theta \,d\theta \wedge d\phi .
\end{eqnarray*}%
\begin{equation*}
\wt R_{\phantom{1}3}^{1}=\frac{1}{2} \wt R_{\phantom{1}313}^{1}\,e^{1}\wedge e^{3}
\end{equation*}%
\begin{equation*}
\wt R_{\phantom{1} 313}^{1}=\frac{4}{f}=\tn{const}. %\cot \theta 
\end{equation*}
%and
%\begin{equation*} \int_{S_6} \wt R_{2323}\,dV = 4\pi %\gamma  ~~~~~~~~~~\tcr{\gamma ?}
%\end{equation*}
In particular we see that the spheres at the two endpoints $r=0 ~\tn{or}~ \pi $ has 
constant positive curvature. These are the zero and infinity section of the Hirzebruch fibration or twisted sphere bundle map. Note that these two sections both have the same diameter since $f(0)=f(\pi)$. The general radius is ${rad(r)}=1/\kappa^2=f^2/16$.

\begin{figure}[h] \bct 
\includegraphics[width=.7\textwidth]{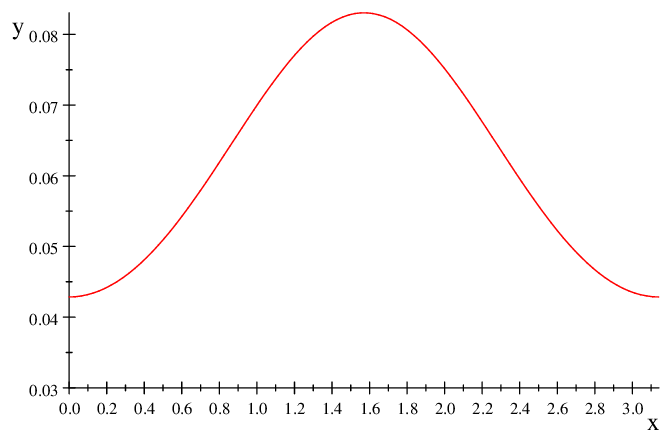}\\ 
\ect
  \caption{\em Radius function $rad(r)$ for $\phi\theta$ spheres. } \label{radius}
\end{figure}

\vspace{.05in}

%$\tcb{\mathbf{S_7:}}$ We take $r=0 \tn{or} \pi $ fixed, and vary $\psi ,\theta $ to 
%obtain an embedded cylinder. If we pass beyond the chart and vary 
%$\theta\in [0,2\pi]$ we get a {\em torus}.  
%Metric on this torus reads 

%\ni 
We can summarize some of the properties we have observed as follows. 

\begin{theorem}
All of the surfaces $S_1 \cdots S_6$ are totally geodesic, hence minimal.
\end{theorem}

%The surfaces $S_5$ and $S_6$ inside the Page space are not totally geodesic. ..... 
%constant curvature - 0 +. 
\begin{theorem} We can list the constant Gaussian curvature surfaces in our list as follows, 
\beg{enumerate}
\item The tori in the classes $S_1$ and $S_3$ are flat. 

\item The spheres in the class $S_6$ are of constant positive curvature. 
\end{enumerate}

\end{theorem}

%\section{Curvature of submanifolds}

%In this section %Now 
%we are going to compute the curvature of % analyze 
%these surfaces.

%\newpage

\section{Three dimensional minimal submanifolds}\label{sechypersurfaces}

In this section we will consider the $3$ dimensional submanifold by fixing a variable 
in our coordinate system.

\vspace{.05in}

$\tcb{\mathbf{N_1:}}$ Easiest case is obtained when we fix the coordinate 
$r=r_0\in (0,\pi)$.  
This way we get totally geodesic $3$-spheres endowed with {\em Berger sphere} metrics. 
At the endpoints $r_0=0,\pi$ this restriction gives round $2$-spheres. 
From (\ref{page1}) 
 we get 
%$$g_{N_1}:=f(r_0)(\sigma _{1}^{2}+\sigma _{2}^{2})+\frac{C\sin ^{2}r_0}{V(r_0)}%
%\sigma _{3}^{2}=\wt\sigma _{1}^{2}+\wt\sigma _{2}^{2}+\frac{C\sin ^{2}r_0}{V(r_0)}%
%\sigma _{3}^{2}$$

{ \renewcommand*{\arraystretch}{2}
$$\begin{array}{rcl}
g_{N_1} &:=& f(r_0)(\sigma _{1}^{2}+\sigma _{2}^{2})+\frac{C\sin ^{2}r_0}{V(r_0)}%
\sigma _{3}^{2} \\
&=&\wt\sigma _{1}^{2}+\wt\sigma _{2}^{2}+\frac{C\sin ^{2}r_0}{V(r_0)}%
\sigma _{3}^{2}.
\end{array}$$
}

$$|\sigma_1|=|\sigma_2|=1/\sqrt{f(r_0)}, ~~~
|\sigma_3|={U \over D \sin r_0}.$$

\ni Then using the orthonormal coframe $\{e^{1},e^{2},e^{3}\}$ %=\{\frac{h}{2}\sin \theta d\phi ,
%\frac{D\sin r}{2U}(d\psi +\cos \theta d\phi ),\frac{h}{2}d\theta \}.
connection 1-forms can be computed from the structure equations as follows,
 $$\hspace{-6mm}\wt\omega^3_{\phantom{1}2}=D U^{-1} h^{-2} \sin r \,e^1,~~\wt\omega^1_{\phantom{1}3}= 2h^{-1}\cot\theta\,e^1 - DU^{-1}h^{-2}\sin r\,e^2,
 ~~\wt\omega^2_{\phantom{1}1}=DU^{-1}h^{-2}\sin r\,e^3$$ 
%We know that the $1$-forms of the ambient connection
%$$w_{3}^{2}=2h^{-1}\cot \theta e^{2},w_{3}^{1}=0,w_{2}^{1}=0,$$ 
which are the same as the one (\ref{ambientconnectiononeforms}) of the ambient space,
%in the ambient space
so that these type of submanifolds are totally geodesic and hence minimal. 
Computing the curvatures, 
\[
\wt R_{\phantom{1} 212}^{1}=0,%
\begin{array}{rr}
& 
\end{array}%
\wt R_{\phantom{1} 313}^{1}=0
\]%
\[
\wt R_{\phantom{1} 323}^{2}=4h^{-2}=const.
\]
So these are constant scalar curvature but non-Einstein (i.e. non-constant sectional curvature in this dimension) spaces. 
\vspace{.05in}

$\tcb{\mathbf{N_2:}}$ Keeping $\theta$ fixed we obtain the Seifert fibered $3$-manifold 
$S^1\times S^2$. To see this note that $\psi$, $\phi$ traces a torus for $0<r<\pi$. At the endpoints the Hopf $\psi$-circles shrink. Looking in the $r,\psi$ direction we have 
an American football (i.e. the spheres of $S_5$), 
on which $\psi$ coordinates draw circles. This is the most direct way to see it. 
Alternatively, start with $I\times T^2$ which one can imagine as a thickened cylinder, 
whose inside is to be glued to its outside. Before gluings lets perform the shrinkings. 
At the endpoints, a set of circles shrink, so the top and 
bottom of the thickened cylinder must be shrank to circles. So this is a thickened 
cylinder but with sharpened ends. If we glue inside out, the neighborhood of these sharpened ends become solid tori. Now we have two solid tori glued along their boundary. 
Solid torus is the product of circle times disk $S^1\times D^2$, and when boundaries 
of two copies are identify, that means identifying the boundary of two disks which has to 
yield a sphere hence sphere times a circle. This solid tori decomposition is actually 
called the Heegard decomposition of this 3-manifold in the $3$-manifold theory literature \cite{rolfsen}. 
A third way of interpretation can be done using the capping process which is described in detail in \cite{handle}. Shrinking a circle is the same operation as attaching a $2$-disk which is called {\em capping}. 
Then it becomes easier to see the two solid tori at the end points. We have a $I\times T^2$, and only the two tori at the to endpoints of this manifold is filled to yield 
solid tori. The rest can be seen as gluing solid tori from a distance. 
%\vspace{.05in} 
The metric on this $3$-manifold is, 

$$\hspace{-10mm}g_{N_2}=Vdr^2+\left\{ {f\over 4}\sin^2\theta +{C\sin^2r \cos^2\theta\over 4V(r)} \right\}\hspace{-1mm}d\phi^2+ {C\sin^2r\over 4V(r)}d\psi^2
+ {C\sin^2r \cos\theta\over 4V(r)}
(d\psi\otimes d\phi+d\phi\otimes d\psi).$$

\ni Then working with the orthonormal coframe 
$\{e^{0},e^{1},e^{2}\}$ %=\{Udr,\frac{h}{2}\sin \theta d\phi ,
%\frac{D\sin r}{2U}(d\psi +\cos \theta d\phi )\}.
we conclude that 
$$\wt\omega^1_{\phantom{1}0}=U^{-1}h^{-1}\dot h \,e^1,~~\wt\omega^2_{\phantom{1}0}=(U^{-1}\ \cot r - U^{-2}\dot U )\,e^2 ~~\tn{and}~~ \wt\omega^1_{\phantom{1}2}=0,$$ 
which is the same as the one in the ambient space
so that this submanifold is totally geodesic. 
%\vspace{.05in} \ni 
Computing the curvature 2-forms,%
\[
\wt R_{\phantom{1} 101}^{0}=2h^{-1}U^{-3}(\ddot hU-\dot h\dot U)
,
\]%
\[
\wt R_{\phantom{1} 202}^{0}=U^{-3}\{2\dot U-U\ddot U-2U\dot U\cot r-U\},%
\begin{array}{rr}
& 
\end{array}%
\wt R_{\phantom{1} 212}^{1}=0.
\]

\vspace{.05in}

$\tcb{\mathbf{N_3:}}$ Keeping $\phi$ fixed we again obtain copies of $S^1\times S^2$ 
similar to the previous case. The metric is, 

$$g_{N_3}=Vdr^2+{C\sin^2r \cos\theta\over 4V(r)}
(d\psi\otimes d\phi+d\phi\otimes d\psi)+{f\over 4}d\theta^2.$$

\ni Then we have the orthonormal coframe 
$$\{e^0,e^2,e^3\}=\{Udr,\frac{D\sin r}{2U}d\psi ,\frac{h}{2}d\theta \},$$
we conclude that 
$$\wt\omega^2_{\phantom{1}0}=(U^{-1}\ \cot r - U^{-2}\dot U )\,e^2,~~ \wt\omega^3_{\phantom{1}0}=0 ~~\tn{and}~~ 
\wt\omega^2_{\phantom{1}3}=0,$$  
which is the same as the one in the ambient space
so that this submanifold is totally geodesic. 
%\vspace{.05in} \ni 
Computing the curvatures, 
\[
\wt R_{\phantom{1} 020}^{2}=U^{-4}\{ \dot UU+\dot U\ \cot r+ 2\dot UU\ \cot r+U^2-3\dot U^2\},
\]%
\[
\wt R_{\phantom{1} 303}^{0}=0,%
\begin{array}{rr}
& 
\end{array}%
\wt R_{\phantom{1} 323}^{2}=0.
\]

\vspace{.05in}

$\tcb{\mathbf{N_4:}}$ Keeping $\psi$ fixed gives us a copy of the {\em cylinder} $S^2\times I$ with the following metric,  
$$g_{N_4}=Vdr^2+\left\{ {f\over 4}\sin^2\theta 
+{C\sin^2r \cos^2\theta\over 4V(r)} \right\}\hspace{-1mm}d\phi^2
+{f\over 4}d\theta^2.$$

\ni Then using the orthonormal coframe 
$\{e^{0},e^{2},e^{3}\}$ %=\{Udr,\frac{D\sin r}{2U}(d\psi +\cos \theta d\phi ),\frac{h}{2}d\theta \}.
we conclude that 
$$\hspace{-14mm}\wt\omega^1_{\phantom{1}0}=U^{-1}h^{-1}\dot h \,e^1 + 2^{-1}U^{-1}h^{-1}\dot h \sin^{-1}\theta\,e^3, ~~  
\wt\omega^3_{\phantom{1}0}=2^{-1}U^{-1}h^{-1}\dot h \sin^{-1}\theta\,e^1, ~~ \wt\omega^3_{\phantom{1}1}=D U^{-1} h^{-2} \sin r \,e^1$$ 
which is the same as the three in the
ambient space so that this submanifold is totally geodesic. Computing the
curvature,%
\[
\wt R_{\phantom{1} 202}^{0}=-2h^{-1}U^{-3}(\ddot hU-\dot h\dot U),
\]%
\[
\wt R_{\phantom{1} 303}^{0}=-h^{-1}U^{-3}(\ddot hU-\dot h\dot U),%
\begin{array}{rr}
& 
\end{array}%
\wt R_{\phantom{1} 323}^{2}=\dot h^{2}h^{-2}U^{-2}.
\]

\vspace{.05in}

%\newpage

\section{Relations with the Yamabe Problem and other metrics on Hirzebruch surfaces}\label{secyamabe}

%The 
We define the {\em normalized Einstein-Hilbert functional} on a smooth 4-manifold by, %is defined as,   

$$E(g)=vol_g^{-1/2} \int_MR_g\,dV_g.$$

\ni Then for a conformal class $[g]=\{fg\,|\, f:M\to\mbbr^+%, f 
~\tn{is}~ C^\infty\}$ of a metric $g$, we define the Yamabe constant of this class by
$Y_{[g]}:=\tn{inf}\{E(e^f g) \,|\, f \,\tn{is}\, C^\infty\}$. This infimum is shown to be achieved by a metric by Yamabe and we call the minimum as the {\em Yamabe minimizer} of the conformal class.  As a result by Aubin, the space of Yamabe constants are bounded from above by the one of round spheres, so that %This function is 
$Y_{[g]}\leq E(S^4_{Round})=24\pi \sqrt{2/3}  \approx 61.562393$. %43.823235$.  
Consequently it has a supremum 
on a manifold $YM(M)=\sup_{g\in \mathcal M}Y_{[g]}$ called the {\em Yamabe Invariant} of the smooth manifold. Here $\mathcal M$ denotes the space of all smooth metrics on the manifold $M$. %$$ \sup_{a\in A} \inf_{b\in B} \mathrm{d}(a,b) $$
We refer to the survey \cite{cleinsteinyamabe} for details.  
The conjectural value of the Yamabe invariant of the Page manifold 
$\mbb{CP}_2\sharp\,\overline{\mbb{CP}}_2$ is 
$12\sqrt{2}\pi\approx 53.314598$ 
which is that of the complex projective plane by \cite{lebrunperturbedcp2}. 
Now, let's calculate the volume of Page metric. Using the coordinates 
(\ref{eulercoordinates}) we have,

\begin{equation*}
g_{ij}=\left( \begin{array}{cccc}
V & 0 & 0 & 0 \\ 
0 & \frac{f}{4}\sin ^{2}\theta +\frac{C\sin ^{2}r\cos ^{2}\theta }{4V} & 
\frac{C\sin ^{2}r\cos \theta }{4V} & 0 \\ 
0 & \frac{C\sin ^{2}r\cos\theta }{4V} & \frac{C\sin ^{2}r}{4V} & 0 \\ 
0 & 0 & 0 & \frac{f}{4}
\end{array}
\right)
\end{equation*}

\ni and then,

%$$\left\vert g_{ij}\right\vert =V\{(\frac{f}{4}\sin ^{2}\theta +\frac{C\sin
%^{2}r\cos ^{2}\theta }{4V})\frac{C\sin ^{2}r}{4V}-(\frac{C\sin
%^{2}r\cos\theta }{4V})^{2}\}{f\over 4}$$

$$\det g_{ij}={C\over 64}f^2\sin^2r\sin^2\theta$$

\ni then the volume form becomes,

$$dV={D\over 8}f(r) \sin r \sin \theta \,dr\wedge d\phi \wedge d\psi\wedge d\theta.$$

%\[dV=\frac{C}{64}f^{2}\sin ^{2}r\sin ^{2}\theta \,dr\wedge d\phi \wedge d\theta\wedge d\psi .\]%

\ni After all taking the variables $\theta\leq\pi$ and $\phi,\psi\leq 2\pi$
 we get the volume,

$$Vol_{g_{page}} = \int_M dV%={16D\pi^2\over 3}{3-a^2\over 3+6a^2-a^4}.
={16\pi^2\over 3(3+a^2)}{3-a^2\over 3+6a^2-a^4}\approx 14.38828. %28.77657.
%\frac{\pi ^{4}}{4}\frac{a^{4}-4a^{2}+16}{(a^{2}+3)^{2}(3+6a^{2}-a^{4})^{2}}
$$
The scalar curvature of the Page metric can be computed to be 
$12(1+a^2)\approx 12.952$. So the Yamabe function takes the following value. 
$$E(g_{Page})=12(1+a^2)vol^{1/2} %\approx 69.481\cdots
\approx 49.13.
%=12(1+a^2)4\sqrt{D}\pi{ \sqrt{3-a^2} \over \sqrt{3+6a^2-a^4} }=120.3441267. 
%6(1+a^2)\pi^2 {\sqrt{a^4-4a^2+16}\over (a^2+3)(3+6a^2-a^4)}\approx 23.694254.
$$  
Since this is an Einstein metric, by \cite{obata71einsteinyamabe} it is a Yamabe minimizer. So that this implies $YM(\mbb{CP}_2\sharp\,\overline{\mbb{CP}}_2)\geq  49.13$. % 69.481$. 
%\vspace{.05in}
On the other hand, the Einstein-Hilbert action of an Otoba metric \cite{otoba} 
on the Page manifold with scalar curvature $R$ is given by
$$E(g_1(R))=2\sqrt[4]{2}\pi R \sqrt{\tn{Arcsin}(k)}, 
~~k^2=(1+\beta/\sqrt{2+\beta^2})/2, 
~~\beta=(8-R)/2.$$
This function indeed takes arbitrary small and large numbers for varying the values of $R$. Since Bach-flat metrics and Otoba metrics are mutually exclusive, 
Page metric is not in the set of Otoba metrics, neither in its conformal rescalings. Due to the theorem of Obata there is no other constant scalar curvature metric in the conformal class of an Einstein metric except for $S^4$. 
Instead, one of the Otoba metrics conformally equivalent to a metric in the K\"{a}hler class of the Page metric.

\newpage

\appendix
\section{Linear Structure Equations for Connection 1-forms}\label{ode}
In this appendix we are going to provide the solution to the 
system of differential equations providing the connection $1$-forms. For 
%this particular equation,
the equation we are dealing with, the simplest particular solution techniques actually do 
 not work. 
So one has to start with the most general form for the possibilities. %solutions. 
That is the reason we decided to include it here. 
We start with writing down the structure equation system, 

\vspace{5mm}   

{ \renewcommand*{\arraystretch}{1.6}  
$\hspace{-2.5cm}\begin{array}{c@{}c@{}cr}
0=&0&+&
~~~0~~~+
\omega^0_{\phantom{1} 1}\wedge e^1 + 
\omega^0_{\phantom{1} 2}\wedge e^2 + 
\omega^0_{\phantom{1} 3}\wedge e^3    ~~~~(0)\\

0=&U^{-1}h^{-1}\dot h\, e^{01}+2h^{-1}\cot\theta\, e^{31}&+&
\omega^1_{\phantom{1} 0}\wedge e^0 + 
~~~0~~~+
\omega^1_{\phantom{1} 2}\wedge e^2 + 
\omega^1_{\phantom{1} 3}\wedge e^3    ~~~~(1)\\

0=&-\sin^{-1}r (U^{-1}\sin r)_r\, e^{20} - 2DU^{-1}h^{-2}\dot h^{-1}\sin r\, e^{31}  &+&
\omega^2_{\phantom{1} 0}\wedge e^0 + 
\omega^2_{\phantom{1} 1}\wedge e^1 + 
~~~0~~~+
\omega^2_{\phantom{1} 3}\wedge e^3    ~~~~(2)\\

0=&U^{-1}h^{-1}\dot h \sin^{-1}\theta\, e^{01} &+&
\omega^3_{\phantom{1} 0}\wedge e^0 + 
\omega^3_{\phantom{1} 1}\wedge e^1 + 
\omega^3_{\phantom{1} 2}\wedge e^2 + ~~~0.~~ ~~~~(3)

\end{array}$
}

\vspace{5mm}   
\ni which involves  a system of {$4$ equations and $6$ unknowns.} 
Using $(0)$ we label the first three forms as follows. 
$$\omega^1_{\phantom{1} 0}=\alpha\,e^1+\beta\,e^2+\gamma\,e^3,
~~\omega^2_{\phantom{1} 0}=\beta\,e^1+\varepsilon\,e^2+\delta\,e^3,
~~\omega^3_{\phantom{1} 0}=\gamma\,e^1+\delta\,e^2+\theta\,e^3.$$

\ni Imposing each of them, the related other equation $(1,2,3)$ we reach,
$$ \hspace{-7mm}
\omega^1_{\phantom{1} 0}=U^{-1}h^{-1}\dot h\,e^1+\beta\,e^2+\gamma\,e^3,
~~\omega^2_{\phantom{1} 0}=\beta\,e^1+\sin^{-1}r(U^{-1}\sin r)_r\,e^2+\delta\,e^3,
~~\omega^3_{\phantom{1} 0}=\gamma\,e^1+\delta\,e^2.$$

\ni Using $(3,2)$ we start with the form 
$\omega^3_{\phantom{1} 2}=-\delta\,e^0+\rho\,e^1$ and 
$\omega^1_{\phantom{1} 2}=\sigma\,e^0+\tau\,e^3$ using $(2)$. 
Next, imposing $(3)$ we start with the following form, 
$$\omega^3_{\phantom{1} 1}=
(\gamma-U^{-1}h^{-1}\dot h\sin^{-1}\theta)\,e^0+\zeta\,e^1+\rho\,e^2.$$

\ni We will work on equation $(1)$. Matching the $03$-terms we can compute 
$\gamma=2^{-1}U^{-1}h^{-1}\dot h\sin^{-1}\theta$. Comparing the $13$-terms 
we get $\zeta=-2h^{-1}\cot\theta$ and comparing $23$-terms we get $\rho=-\tau$.   
We continue comparisons within the same equation by $02$-terms to get $\sigma=\beta$. 
Comparison within the next equation $(2)$ yields $\beta=0$ and 
$\tau=-DU^{-1}h^{-2}\sin r$. Collecting these we able to settle the forms, 
$$\omega^3_{\phantom{1} 1}=-2^{-1}U^{-1}h^{-1}\dot h\sin^{-1}\theta\,e^0+
-2h^{-1}\cot\theta\,e^1+DU^{-1}h^{-2}\sin r\,e^2,$$ %~~~~
$$\omega^1_{\phantom{1} 2}=-DU^{-1}h^{-2}\sin r\,e^3.$$
To figure out $\delta$ we consult to the final equation $(3)$. 
Comparing the $02$-terms here yields $\delta=0$. 
%An alternative for the middle couple in (\ref{vierbein}) 
%yielding more cumbersone computations is, 

%$$\hspace{-10mm}\{ {D\sin r\over 2U}d\psi,
%{D\sin r\over 2U\sqrt{1+C^{-1}Vf\cos^{-2}\theta\sin^2\theta\sin^{-2}r} } 
%(d\psi+(C^{-1}Vf\cos^{-1}\theta\sin^2\theta\sin^{-2}r+\cos\theta)d\phi)\}.$$ 
%$$\hspace{-20mm}\{ {Dh\over 2\sqrt{C \cot^2\theta+Vf\csc^2r}} d\psi,
%{D\sin r\over 2U\sqrt{1+C^{-1}Vf\tan^2\theta\sin^{-2}r} } 
%(d\psi+(C^{-1}Vf\tan\theta\sin\theta\sin^{-2}r+\cos\theta)d\phi)\}.$$ 

%which yields more cumbersone computations.

\vspace{1cm}

%\newpage

\bibliography{page}{}
\bibliographystyle{alphaurl}

\vspace{.05in}

%\vspace{2cm}

%{\small
%\beg{flushleft}
%\textsc{%Mathematics Department,
%Michigan State University,
%East Lansing, MI 48824, USA}\\
%\textit{E-mail address :} \texttt{\textbf{kalafat@math.msu.edu}}
%\end{flushleft}
%}

{\small
\beg{flushleft}
\textsc{Orta mh. Z\"ubeyde Han\i m cd. No 5-3 Merkez 74100 Bart\i n, T\" urk\'{i}ye.}\\
\textit{E-mail address:} \texttt{\textbf{kalafat@\,math.msu.edu}}
\end{flushleft}
}

{\small \beg{flushleft}
\textsc{Amasya \" Un\' ivers\'ites\'i, Merz\' ifon Myo, 05300, T\"urk\'iye.}\\
\textit{E-mail address:} \texttt{\textbf{ramazan.sari@\,amasya.edu.tr}}

%Amasya University, Merzifon Vocational Schools, 05300, Amasya, Turkey

\end{flushleft}}

\end{document}